\documentclass[12pt]{amsart}
\usepackage{amssymb}

\usepackage{graphicx}

\usepackage[]{amsfonts}

\def\underset#1#2{\mathrel{\mathop{\kern0pt #2}\limits_{#1}}}

\def\overset#1#2{\mathrel{\mathop{\kern0pt #2}\limits^{#1}}}

\def\couleur(#1 #2 #3)
	{
\special{ps::[end]}}

\def\sqr#1#2{{\vcenter{\vbox{\hrule height.#2pt
             \hbox{\vrule width.#2pt height#1pt \kern#1pt
             \vrule width.#2pt}
             \hrule height.#2pt}}}}

\def\st{\mathinner{\mkern1mu\raise1pt\hbox{.}				
		   \mkern1mu\raise4pt\hbox{.}
		   \mkern1mu\raise1pt\hbox{.}
		 }
         }

\def\bx#1{\setbox1=\hbox{\kern3pt{#1}\kern3pt}				
 \dimen1=\ht1 \advance\dimen1 by 3pt \dimen2=\dp1 \advance\dimen2 by 3pt
 \setbox1=\hbox{\vrule height\dimen1 depth\dimen2\box1\vrule}%
 \setbox1=\vbox{\hrule\box1\hrule}%
 \advance\dimen1 by .4pt \ht1=\dimen1
 \advance\dimen2 by .4pt \dp1=\dimen2 \box1\relax}

\def\k#1{\kern#1em}
\def\vci{\vrule  width.02em height1.47ex depth-.0ex}				
\def\11{{\rm\k{.2}\vci\k{-.37}1}}

\newtheorem{Theorem}{Theorem}[section]
\newtheorem{Proposition}[Theorem]{Proposition}
\newtheorem{Lemma}[Theorem]{Lemma}

\begin{document}
\title{Numerical Solution of a nonlinear reaction-diffusion problem in 
the case of HS-Regime}
\author{Marie-Noelle Le Roux}
\address{IMB, Institut de Math\'ematiques de Bordeaux, 
UMR5251\\Universit{\'e} Bordeaux1, 33405 Talence Cedex, FRANCE}
\email{Marie-Noelle.Leroux@math.u-bordeaux1.fr}
\maketitle
\begin{abstract} {
In this paper, the author propose a numerical method to compute the 
solution of the Cauchy problem: $w_{t}-(w^{m}w_{x})_{x}=w^{p}$, the initial 
condition is a nonnegative function with compact support, $m>0,\ 1<p<m+1$. 
The problem is split in two parts: A hyperbolic term solved by using 
the Hopf and Lax formula and a parabolic term solved by a backward linearized 
Euler method in time and a finite element method in space. Estimates 
of the numerical solution are obtained and it is proved that any numerical 
solution is unbounded.\ \par
}\end{abstract}
\ \par
\ \par
\ \par
\section{Introduction}
\setcounter{equation}{0}\ \par
In this paper, we study a numerical method to compute the solution of 
the Cauchy problem  :\ \par

\begin{equation} 
\begin{matrix}{w_{t}-(w^{m}w_{x})_{x}=w^{p},\ t>0,\ x\in {\mathbb{R}}}\cr 
{w(x,0)=w_{0}(x)\geq 0,\ x\in {\mathbb{R}}}\cr \end{matrix}\label{HSregime0}
\end{equation} \ \par
$w_{0}$ is a function with compact support, $m>0$ , $1<p<m+1$.\ \par
Samarskii et al ~\cite{Samarskii1}, see also~\cite{Bandle}, ~\cite{Bebernes}, 
 ~\cite{Galaktionov1}, ~\cite{Galaktionov2}, ~\cite{Galaktionov3}, ~\cite{Galaktionov4} 
, ~\cite{Sacks}  have obtained theoretical results on this problem. In 
the case $1<p<m+1$, the numerical solution blows up in finite time and 
there is no localization (HS-regime) that is $u(t,x){\longrightarrow}\infty \ in\ {\mathbb{R}}\ if\ t{\longrightarrow}T_{0}$.\ 
\par
\ \par
{\hskip 1.2em}A numerical method to solve ~(\ref{HSregime0}) has been 
proposed in the case of S-regime ($\displaystyle p=m+1$) in ~\cite{mnLeRoux6} 
and in the case of LS-regime ($p>m+1)$  ~\cite{mnLeRoux7}. If we denote 
$\omega _{L}=\left\lbrace{\left.{x\in {\mathbb{R}}/u(T_{0}^{-},x)=\infty 
}\right\rbrace }\right. $, in the case of S-regime, $L^{*}=mes(\omega _{L})$ 
is positive while in the case of LS-regime , $L^{*}=0$.  The problem 
is solved by using a splitting method; for that, it is more convenient 
to work with the function $u=w^{m}$. Problem ~(\ref{HSregime0}) may be 
written:\ \par

\begin{equation} 
\begin{matrix}{u_{t}-\frac{1}{m}u_{x}^{2}-uu_{xx}=mu^{q+1},\ t>0,\ x\in 
{\mathbb{R}}}\cr {u(x,0)=u_{0}(x)=w_{0}^{m}(x),\ x\in {\mathbb{R}}}\cr 
\end{matrix}\label{HSregime1}
\end{equation}   \ \par
with $q=\frac{p-1}{m},\ m>0,\ q<1$.\ \par
This problem is split in two parts: a hyperbolic problem which will be 
solved exactly at the nodes at each time step and allows the extension 
of the domain and a parabolic problem which will be solved by a backward 
linearized Euler method which allows the blow up of the solution. \ \par
In ~\cite{mnLeRoux6},~\cite{mnLeRoux7}, the convergence of the scheme 
has been proved in the cases $q=1$ and $q>1$. It has also been proved 
that for $q=1$, the numerical solution blows up in  finite time for any 
initial condition and that its support remains bounded if the initial 
condition is smaller than a self-similar solution. In the case, $1\leq q<\frac{m+2}{m},$ 
any numerical solution is unbounded while for $q>\frac{m+2}{m}$  ($p>m+3$) 
 if the initial condition is sufficiently small, a global solution exists 
and if $q\geq \frac{m+2}{m}\ $ for large initial condition, the solution 
blows up in a finite time. We  observe numerically that in any case, 
the  unbounded solution is strictly localized and blows up in one point 
and that for $q=\frac{m+2}{m\ }\ $also, if the initial condition is sufficiently 
small, a global solution exists. \ \par
Here, we generalize this method to the case $q<1.$\ \par
\ \par
An outline of the paper is as follows:\ \par
{\hskip 1.8em}In Section 2, we present the numerical scheme. In Section 
3, we obtain estimates of the approximate solution and of its derivative 
in space. In Section 4, we prove that if $q<1$, ($p<m+1$)) any numerical 
solution is unbounded  \ \par
\ \par
\section{Definition of the numerical solution}
\setcounter{equation}{0}\ \par
In order to solve problem ~(\ref{HSregime1}), we separate it in two parts: 
 a hyperbolic  problem\ \par

\begin{equation} 
u_{t}-\frac{1}{m}(u_{x})^{2}=0,\ x\in {\mathbb{R}},\ t>0\label{HSregime2}
\end{equation} \ \par
and a parabolic problem\ \par

\begin{equation} 
u_{t}-uu_{xx}=mu^{q+1},\ x\in {\mathbb{R}},\ t>0\label{HSregime3}
\end{equation} \ \par
\ \par
{\hskip 1.8em}We denote by $\Delta t_{n}$ the time increment between 
the time levels $t_{n}$ and $t_{n+1},\ n\geq 0$ and by $u_{h}^{n}$  the 
approximate solution at the time level $t_{n}$. This solution will be 
in a finite-dimensional space which will be defined below.\ \par
{\hskip 1.8em}Without loss of generality, we can assume that the initial 
condition is a continous function with a symmetric compact support $[-s_{0},s_{0}]$. 
 Let$\ N\in {\mathbb{N}}$,  the space step $h$  is defined by $h=\frac{s_{0}}{N};$ 
we note $x_{i}=ih,\ i\in {\mathbb{Z}},\ I_{i}=(x_{i-1},x_{i})$ and we 
define the finite dimensional  space $V_{h}^{0}$ by\ \par
\ \par

\begin{equation} 
V_{h}^{0}=\left\lbrace{\left.{\phi _{h}\in C^{0}({\mathbb{R}})/\phi _{h}(x)=0,\ 
x\not\in ]-s_{0},+s_{0}[,\ \phi _{h\mid I_{i}}\in P_{1},\ i=-N+1,N}\right\rbrace 
}\right. 
\end{equation} \ \par
For $\phi _{h}\in V_{h}^{0}$, we note $\phi _{i}=\phi _{h}(x_{i}),\ i\in {\mathbb{Z}}$ 
and for any function $v\in C^{0}({\mathbb{R}})$ with compact support 
$[-s_{0},+s_{0}]$, we define its interpolate by $\pi _{h}^{0}v\in V_{h}^{0}$ 
and $\pi _{h}^{0}v(x_{i})=v(x_{i}),\ i\in {\mathbb{Z}}$.\ \par
{\hskip 1.8em}The support of the solution $u_{h}^{n}$ will be denoted 
$[-s_{-}^{n},\ s_{+}^{n}]$ and will be computed at each time level by 
solving ~(\ref{HSregime2}).\ \par
{\hskip 2.7em}We denote\ \par

\begin{displaymath} 
N_{-}^{n}=\left[{\frac{s_{-}^{n}}{h}}\right] ,\ N_{+}^{n}=\left[{\frac{s_{+}^{n}}{h}}\right] 
,\ h_{-}^{n}=s_{-}^{n}-(N_{-}^{n}-1)h,\ h_{+}^{n}=s_{+}^{n}-(N_{+}^{n}-1)h\end{displaymath} 
\ \par
(for $x\in {\mathbb{R}}$, $\left[{x}\right] $ denotes the greatest integer 
less than $x$), so we get: $h\leq h_{-}^{n},h_{+}^{n}<2h$.\ \par
\ \par
We then define the finite-dimensional space $V_{h}^{n}$ by:\ \par

\begin{displaymath} 
V_{h}^{n}=\left\lbrace{\left.{\begin{matrix}{\phi _{h}\in C^{0}({\mathbb{R}})/\phi 
_{h}(x)=0,\ x\not\in ]-s_{-}^{n},s_{+}^{n}[,}\cr {\ \phi _{h\mid I_{i}}\in 
P_{1},i=-N_{-}^{n}+2,N_{+}^{n}-1}\cr {\phi _{h\mid (-s_{-}^{n},x_{-N_{-}^{n}+1})},\phi 
_{h\mid (x_{N_{+}^{n}-1},s_{+}^{n})}\in P_{1}}\cr \end{matrix}}\right\rbrace 
}\right. \end{displaymath} \ \par
\ \par
and denote by $\pi _{h}^{n}v$ the Lagrange interpolate in $V_{h}^{n}$ 
of a function $v\in C^{0}({\mathbb{R}})$ with a compact support in $[-s_{-}^{n},\ s_{+}^{n}]$.\ 
\par
{\hskip 1.8em}The solution $u_{h}^{n+1}$ at the time level $t_{n+1}\ $is 
computed in two steps: knowing $u_{h}^{n}$, we compute first an approximate 
solution of ~(\ref{HSregime2}) which will be denoted $u_{h}^{n+\frac{1}{2}}$ 
; then starting  with this intermediate value, we compute an approximate 
solution of ~(\ref{HSregime3}) at time level $t_{n+1}$.\ \par
{\hskip 1.5em}Then if $S$ is the semi-group operator associated with 
~(\ref{HSregime2}), we define 
\begin{displaymath} 
u^{n+\frac{1}{2}}=S(\Delta t_{n})u_{h}^{n}.\end{displaymath}   The support 
of this function will be denoted $[-s_{-}^{n+1},s_{+}^{n+1}]\ $ and the 
interpolate of $u^{n+\frac{1}{2}}$ in $V_{h}^{n+1}$ will be denoted $u_{h}^{n+\frac{1}{2}}=\pi _{h}^{n+1}u^{n+\frac{1}{2}}$.\ 
\par
{\hskip 1.8em}Then starting with $u_{h}^{n+\frac{1}{2}}$, the function 
$u_{h}^{n+1}$ is obtained by solving ~(\ref{HSregime3})  with a backward 
linearized Euler method in time and a $P_{1}$-finite element method with 
numerical integration in space. This function has the same support as 
$u^{n+\frac{1}{2}}_{h}$.\ \par
\ \par
\subsection{Computation of the solution of the hyperbolic problem}
\ \par
The hyperbolic problem is independent of $q$; we use the same method 
as in~\cite{mnLeRoux6}  for the case $q=1$.It is not necessary in this 
case to use $P_{2}$-interpolation on the last interval since there no 
localization. We use the Hopf and Lax formula which gives explicitely 
the solution to ~(\ref{HSregime2}) with the starting data $u_{h}^{n}\ $at 
the time level $t=t_{n}$. Here, we simply recall the results obtained 
in ~\cite{mnLeRoux6}\ \par
{\hskip 1.5em}\ \par
{\hskip 1.8em}We define the piecewise constant function$v_{h}^{n}$ by\ 
\par

\begin{displaymath} 
v_{i}^{n}=\frac{u_{i}^{n}-u_{i-1}^{n}}{h}\ on\ I_{i}\ ,\ -N_{-}^{n}+2\leq 
i\leq N_{+}^{n}-1\end{displaymath} \ \par

\begin{displaymath} 
v_{N_{+}^{n}}^{n}=-\frac{u^{n}_{N_{+}^{n}-1}}{h_{+}^{n}}\ on\ (x_{N_{+}^{n}-1},s_{+}^{n})\end{displaymath} 
\ \par

\begin{displaymath} 
v^{n}_{-N_{-}^{n}+1}=\frac{u^{n}_{-N_{-}^{n}+1}}{h_{-}^{n}}\ on\ (-s_{-}^{n},\ 
x_{-N_{-}^{n}+1}).\end{displaymath} \ \par
\ \par
Let us denote $r_{n}=\frac{\Delta t_{n}}{h},\ \left\Vert{v}\right\Vert _{s}=\left\Vert{v}\right\Vert 
_{L^{s}({\mathbb{R}})},\ s>0$.\ \par
\begin{Proposition} If the following stability condition \ \par

\begin{equation} 
r_{n}\left\Vert{v_{h}^{n}}\right\Vert _{\infty }\leq \frac{m}{2}
\end{equation}  is satisfied, then the solution $u_{h}^{n+1}$ of ~(\ref{HSregime2}) 
is defined by\ \par

\begin{displaymath} 
u_{i}^{n+\frac{1}{2}}=u_{i}^{n}+\frac{\Delta t_{n}}{m}\left({{\mathrm{max}}(0,-v_{i}^{n},v_{i+1}^{n})}\right) 
^{2},\ -N_{-}^{n}+1\leq i\leq N_{+}^{n}-1\end{displaymath} \ \par

\begin{displaymath} 
s_{+}^{n+1}=s_{+}^{n}-\frac{\Delta t_{n}}{m}v^{n}_{N_{+}^{n}}\end{displaymath} 
\ \par

\begin{displaymath} 
s_{-}^{n+1}=s_{-}^{n}+\frac{\Delta t_{n}}{m}v^{n}_{-N_{-}^{n}+1}\end{displaymath} 
If $N_{+}^{n+1}=N_{+}^{n}+1$, we get \ \par

\begin{displaymath} 
u^{n+\frac{1}{2}}_{N_{+}^{n}}=\left({1-\frac{h}{h_{+}^{n}}}\right) u^{n}_{N_{+}^{n}-1}+\frac{\Delta 
t_{n}}{m}\left({v^{n}_{N_{+}^{n}}}\right) ^{2}\end{displaymath} \ \par
and we get analogous formula at the other end of the support.\ \par
\end{Proposition}
\ \par
\subsection{Computation of the parabolic problem}
\ \par
The approximate solution at $t_{n+1}$ is now obtained by solving problem 
~(\ref{HSregime3}). \ \par
We introduce the approximate scalar product on $V_{h}^{n+1}:\ \forall \phi _{h},\psi _{h}\in V_{h}^{n+1},$\ 
\par

\begin{displaymath} 
(\phi _{h},\psi _{h})_{h}=\frac{1}{2}(h_{-}^{n+1}+h)\phi _{-N_{-}^{n+1}+1}\psi 
_{-N_{-}^{n+1}+1}+h\sum_{i=-N_{-}^{n+1}+2}^{i=N_{+}^{n+1}-2}{\phi _{i}\psi 
_{i}}+\frac{1}{2}(h_{+}^{n+1}+h)\phi _{N_{+}^{n+1}-1}\psi _{N_{+}^{n+1}-1}.\end{displaymath} 
\ \par
\ \par
We define  $u_{h}^{n+1}$ as the solution of the following problem:\ \par

\begin{equation} 
\forall \phi _{h}\in V_{h}^{n+1},\left\lbrace{\begin{matrix}{(u_{h}^{n+1},\phi 
_{h})_{h}+\Delta t_{n}((u_{h}^{n+1})_{x},(u_{h}^{n+\frac{1}{2}}\phi _{h})_{x}=}\cr 
{(u_{h}^{n+\frac{1}{2}},\phi _{h})_{h}+mq\Delta t_{n}\left({\pi _{h}^{n+1}\left({\left({u_{h}^{n+\frac{1}{2}}}\right) 
^{q}u_{h}^{n+1}}\right) ,\phi _{h}}\right) _{h}}\cr {+m(1-q)\Delta t_{n}\left({\pi 
_{h}^{n+\frac{1}{2}}\left({\left({u_{h}^{n+\frac{1}{2}}}\right) ^{q+1}}\right) 
,\phi _{h}}\right) _{h}}\cr \end{matrix}}\right. \label{HSregime5}
\end{equation} \ \par
\ \par
The second member of ~(\ref{HSregime5}) is splitted in two parts in order 
to obtain the $L^{\infty }$-estimate of $u_{h}^{n+1}.$\ \par
This equation may be written:\ \par

\begin{displaymath} 
\left({1-mq\Delta t_{n}\left({u_{i}^{n+\frac{1}{2}}}\right) ^{q}}\right) 
u_{i}^{n+1}+\frac{\Delta t_{n}}{h^{2}}u_{i}^{n+\frac{1}{2}}\left({2u_{i}^{n+1}-u_{i-1}^{n+1}-u_{i+1}^{n+1}}\right) 
=\end{displaymath} \ \par

\begin{displaymath} 
u_{i}^{n+\frac{1}{2}}+m(1-q)\Delta t_{n}\left({u_{i}^{n+\frac{1}{2}}}\right) 
^{q+1},\ -N_{-}^{n+1}+2\leq i\leq N_{+}^{n+1}-2\label{HSregime4}\end{displaymath} 
\ \par
\ \par

\begin{displaymath} 
\left({1-mq\Delta t_{n}\left({u_{N_{+}^{n+1}-1}^{n+\frac{1}{2}}}\right) 
^{q}}\right) u^{n+1}_{N_{+}^{n+1}-1}+\end{displaymath} \ \par

\begin{displaymath} 
\frac{\Delta t_{n}}{h}u^{n+\frac{1}{2}}_{N_{+}^{n+1}-1}\left({\frac{2}{h_{+}^{n+1}}u_{N_{+}^{n+1}-1}^{n+1}-\frac{2}{h+h_{+}^{n+1}}u_{N_{+}^{n+1}-2}}\right) 
\end{displaymath} \ \par

\begin{displaymath} 
=u_{N_{+}^{n+1}-1}^{n+\frac{1}{2}}+m(1-q)\Delta t_{n}\left({u_{N_{+}^{n+1}-1}^{n+\frac{1}{2}}}\right) 
^{q+1}\label{HSregime29}\end{displaymath} \ \par
\ \par

\begin{displaymath} 
\left({1-mq\Delta t_{n}\left({u_{-N_{-}^{n+1}+1}^{n+\frac{1}{2}}}\right) 
^{q}}\right) u^{n+1}_{-N_{-}^{n+1}+1}+\frac{\Delta t_{n}}{h}u^{n+\frac{1}{2}}_{-N_{-}^{n+1}+1}\left({\frac{2}{h_{-}^{n+1}}u^{n+1}_{-N_{-}^{n+1}+1}-\frac{2}{h_{-}^{n+1}+h}u^{n+1}_{-N_{-}^{n+1}+2}}\right) 
\end{displaymath} \ \par

\begin{displaymath} 
=u^{n+\frac{1}{2}}_{-N_{-}^{n+1}+1}+m(1-q)\Delta t_{n}\left({u^{n+\frac{1}{2}}_{-N_{-}^{n+1}+1}}\right) 
^{q+1}\label{HSregime210}\end{displaymath} \ \par
\ \par
\ \par
We get immediately the result:\ \par
\ \par
\begin{Proposition}  \label{HSregime5}If the hypotheses of proposition 
~\ref{HSregime3} are satisfied and if $u_{h}^{n+\frac{1}{2}}$ satisfies\ 
\par

\begin{equation} 
\ mq\Delta t_{n}\left\Vert{u_{h}^{n+\frac{1}{2}}}\right\Vert _{\infty 
}^{q}<1
\end{equation} \ \par
then the solution $\displaystyle u_{h}^{n+1}$ of ~(\ref{HSregime5}) is 
unique and nonnegative.  \ \par
\end{Proposition}
\ \par
\subsection{Properties of the scheme}
\ \par
\begin{Lemma} \label{HSregime8}If the hypotheses of proposition ~\ref{HSregime5} 
are satisfied, then the following estimate holds:\ \par

\begin{equation} 
\left\Vert{u_{h}^{n+1}}\right\Vert _{\infty }\leq \frac{\left\Vert{u_{h}^{n}}\right\Vert 
_{\infty }\left({1+m(1-q)\Delta t_{n}\left\Vert{u_{h}^{n}}\right\Vert 
_{\infty }^{q}}\right) }{1-mq\Delta t_{n}\left\Vert{u_{h}^{n}}\right\Vert 
_{\infty }^{q}}\label{HSregime7}
\end{equation} \ \par
\end{Lemma}
\ \par
Proof: We get immediately from the Hopf and Lax formula that $\left\Vert{u^{n+\frac{1}{2}}_{h}}\right\Vert _{\infty }\leq \left\Vert{u_{h}^{n}}\right\Vert 
_{\infty }$. If we denote $i_{0}$ the index such that $\displaystyle 
u_{\displaystyle i_{0}}^{n+1}=\displaystyle \left\Vert{\displaystyle u_{h}^{n+1}}\right\Vert 
_{\infty }$, we get from ~(\ref{HSregime4}),~(\ref{HSregime29}), ~(\ref{HSregime210}) 
 that\ \par

\begin{displaymath} 
u_{i_{0}}^{n+1}\leq u_{i_{0}}^{n+\frac{1}{2}}\frac{1+m(1-q)\Delta t_{n}\left({u_{i_{0}}^{n+\frac{1}{2}}}\right) 
^{q}}{1-mq\Delta t_{n}\left\Vert{u_{h}^{n}}\right\Vert _{\infty }^{q}}\end{displaymath} 
\ \par
\ \par
which proves the lemma.\ \par
\ \par
We deduce the following theorem:\ \par
\ \par
\begin{Theorem} Under the hypotheses of  proposition ~\ref{HSregime5}, 
the numerical solution exists at least until the time \ \par

\begin{equation} 
T_{1}=\frac{1}{mq\left\Vert{u_{h}^{0}}\right\Vert _{\infty }^{q}}
\end{equation}  \ \par
and the following estimate holds:\ \par

\begin{equation} 
\left\Vert{u_{h}^{n}}\right\Vert _{\infty }\leq \frac{\left\Vert{u_{h}^{0}}\right\Vert 
_{\infty }}{\left({1-mqt_{n}\left\Vert{u_{h}^{0}}\right\Vert _{\infty 
}^{q}}\right) ^{\frac{1}{q}}}\label{HSregime8}
\end{equation} \ \par
\end{Theorem}
\ \par
Proof: This result is proved recurently. It is true for $n=0$. If we 
suppose that we have estimate ~(\ref{HSregime8}) at the time level $t_{n}$ 
, we get from ~(\ref{HSregime7}) ,  at the time \ \par
$t_{n+1:}$\ \par

\begin{displaymath} 
\left\Vert{u_{h}^{n+1}}\right\Vert _{\infty }\leq \left\Vert{u_{h}^{n}}\right\Vert 
_{\infty }\frac{1+m(1-q)\Delta t_{n}\left\Vert{u_{h}^{n}}\right\Vert 
_{\infty }^{q}}{1-mq\Delta t_{n}\left\Vert{u_{h}^{n}}\right\Vert _{\infty 
}^{q}}\end{displaymath} \ \par
or 
\begin{displaymath} 
\left\Vert{u_{h}^{n+1}}\right\Vert _{\infty }\leq \left\Vert{u_{h}^{0}}\right\Vert 
_{\infty }\frac{1-mqt_{n}\left\Vert{u_{h}^{0}}\right\Vert _{\infty }^{q}+m(1-q)\Delta 
t_{n}\left\Vert{u_{h}^{0}}\right\Vert _{\infty }^{q}}{\left({1-mqt_{n+1}\left\Vert{u_{h}^{0}}\right\Vert 
_{\infty }^{q}-}\right) \left({1-mqt_{n}\left\Vert{u_{h}^{0}}\right\Vert 
_{\infty }^{q}}\right) ^{\frac{1}{q}}}\end{displaymath} \ \par
The inequality ~(\ref{HSregime8}) will be satisfied at the time $t_{n+1}\ $if 
:\ \par

\begin{displaymath} 
\left({1-mqt_{n+1}\left\Vert{u_{h}^{0}}\right\Vert _{\infty }^{q}+m\Delta 
t_{n}\left\Vert{u_{h}^{0}}\right\Vert _{\infty }^{q}}\right) \end{displaymath} 
\ \par

\begin{displaymath} 
\leq \left({1-mqt_{n}\left\Vert{u_{h}^{0}}\right\Vert _{\infty }^{q}}\right) 
^{\frac{1}{q}}\left({1-mqt_{n+1}\left\Vert{u_{h}^{0}}\right\Vert _{\infty 
}^{q}}\right) ^{1-\frac{1}{q}}\label{HSregime9}\end{displaymath} \ \par
By using the Taylor formula, we get:\ \par

\begin{displaymath} 
\left({1-mqt_{n}\left\Vert{u_{h}^{0}}\right\Vert _{\infty }^{q}}\right) 
^{\frac{1}{q}}-\left({1-mqt_{n+1}\left\Vert{u_{h}^{0}}\right\Vert _{\infty 
}^{q}}\right) ^{\frac{1}{q}}\geq m\Delta t_{n}\left\Vert{u_{h}^{0}}\right\Vert 
_{\infty }^{q}\left({1-mqt_{n+1}\left\Vert{u_{h}^{0}}\right\Vert _{\infty 
}^{q}}\right) ^{\frac{1}{q}-1}\end{displaymath} \ \par
\ \par
and the inequality~(\ref{HSregime9}) is  satisfied: \ \par
\ \par
\ \par
\begin{Lemma} Under the hypothesis of Proposition(~\ref{HSregime3}), 
we have the estimate:\ \par

\begin{displaymath} 
\left\Vert{v_{h}^{n}}\right\Vert _{\infty }\leq C\end{displaymath} \ 
\par
for $t_{n}\leq T<T_{1}\ $where $C\ $is a constant depending on $T$ and 
$u_{0}$.\ \par
\end{Lemma}
\ \par
\textbf{Proof}: We have the inequality(~\cite{mnLeRoux6} ):$\left\Vert{v_{h}^{n+\frac{1}{2}}}\right\Vert _{\infty }\leq \left\Vert{v_{h}^{n}}\right\Vert 
_{\infty }$.\ \par
It remains to estimate $\left\Vert{v_{h}^{n+1}}\right\Vert _{\infty }$ 
.\ \par
From ~(\ref{HSregime4}), we deduce the following equation satisfied by 
$v_{h}^{n+1}$:\ \par

\begin{displaymath} 
\left({1-mq\Delta t_{n}\left({u_{i}^{n+\frac{1}{2}}}\right) ^{q}}\right) 
v_{i}^{n+1}+\frac{\Delta t_{n}}{h^{2}}\left({v_{i}^{n+1}\left({u_{i}^{n+\frac{1}{2}}+u_{i-1}^{n+\frac{1}{2}}}\right) 
-v_{i-1}^{n+1}u_{i-1}^{n+\frac{1}{2}}-v_{i+1}^{n+1}u_{i}^{^{n+\frac{1}{2}}}}\right) 
\end{displaymath} \ \par

\begin{displaymath} 
=v_{i}^{n+\frac{1}{2}}+mq\frac{\Delta t_{n}}{h}\left({\left({u_{i}^{n+\frac{1}{2}}}\right) 
^{q}-\left({u_{i-1}^{n+\frac{1}{2}}}\right) ^{q}}\right) u_{i-1}^{n+1}+m(1-q)\frac{\Delta 
t_{n}}{h}\left({\left({u_{i}^{n+\frac{1}{2}}}\right) ^{q+1}-\left({u_{i-1}^{n+\frac{1}{2}}}\right) 
^{q+1}}\right) \end{displaymath} \ \par
\ \par
$\displaystyle -N_{-}^{n+1}+3\leq i\leq N_{+}^{n+1}-2$\ \par
\ \par
and we have analogous inequalities for $\displaystyle i=N_{+}^{n+1}-1,\ N_{+}^{n+1},-N_{-}^{n+1}+2,-N_{-}^{n+1}+1$. 
\ \par
By using  ~(\ref{HSregime4}) for $i-1$, we can replace $u_{i-1}^{n+1}\ $in 
the second member by its expression in function of the values of $u_{h}^{n+\frac{1}{2}}$ 
and $v_{h}^{n+1}$ and we get:\ \par

\begin{displaymath} 
\left({1-mq\Delta t_{n}\left({u_{i}^{n+\frac{1}{2}}}\right) ^{q}}\right) 
v_{i}^{n+1}+\end{displaymath} \ \par

\begin{displaymath} 
\frac{\Delta t_{n}}{h^{2}}\left({\left({v_{i}^{n+1}-v_{i+1}^{n+1}}\right) 
u_{i}^{n+\frac{1}{2}}+\left({v_{i}^{n+1}-v_{i-1}^{n+1}}\right) u_{i}^{n+\frac{1}{2}}\frac{1-mq\Delta 
t_{n}\left({\left({u_{i}^{n+\frac{1}{2}}}\right) ^{q}-\left({u_{i-1}^{n+\frac{1}{2}}}\right) 
^{q}}\right) }{1-mq\Delta t_{n}\left({u_{i}^{n+\frac{1}{2}}}\right) ^{q}}}\right) 
\end{displaymath} \ \par

\begin{displaymath} 
=v_{i}^{n+\frac{1}{2}}+mq\frac{\Delta t_{n}}{h}\left({\left({u_{i}^{n+\frac{1}{2}}}\right) 
^{q}-\left({u_{i-1}^{n+\frac{1}{2}}}\right) ^{q}}\right) u_{i}^{n+\frac{1}{2}}\frac{1+m(1-q)\Delta 
t_{n}\left({u_{i-1}^{n+\frac{1}{2}}}\right) ^{q}}{1-mq\Delta t_{n}\left({u_{i-1}^{n+\frac{1}{2}}}\right) 
^{q}}\end{displaymath} \ \par

\begin{displaymath} 
+m(1-q)\frac{\Delta t_{n}}{h}\left({\left({u_{i}^{n+\frac{1}{2}}}\right) 
^{q+1}-\left({u_{i-1}^{n+\frac{1}{2}}}\right) ^{q+1}}\right) \end{displaymath} 
\ \par
\ \par
Let $i_{0}$ the index such that $\displaystyle \displaystyle \left\vert{\displaystyle v_{\displaystyle i_{0}}^{n+1}}\right\vert 
=\displaystyle \left\Vert{\displaystyle v_{h}^{n+1}}\right\Vert _{\infty 
}.\ $From the previous equality, we get:\ \par

\begin{displaymath} 
\left({1-mq\Delta t_{n}\left\Vert{u_{h}^{n}}\right\Vert _{\infty }^{q}}\right) 
\left\vert{v_{i_{0}}^{n+1}}\right\vert \leq \left\vert{v_{i_{0}}^{n+\frac{1}{2}}}\right\vert 
+m(1-q)\frac{\Delta t_{n}}{h}\left\vert{\left({u_{i_{0}}^{n+\frac{1}{2}}}\right) 
^{q+1}-\left({u_{i_{0}-1}^{n+\frac{1}{2}}}\right) ^{q+1}}\right\vert 
\end{displaymath} \ \par

\begin{displaymath} 
+mq\frac{\Delta t_{n}}{h}\left\vert{\left({u_{i_{0}}^{n+\frac{1}{2}}}\right) 
^{q}-\left({u_{i_{0}}^{n+\frac{1}{2}}}\right) ^{q}}\right\vert u_{i_{0}}^{n+\frac{1}{2}}\frac{1+m(1-q)\Delta 
t_{n}\left\Vert{u_{h}^{n+1}}\right\Vert _{\infty }}{1-mq\Delta t_{n}\left\Vert{u_{h}^{n}}\right\Vert 
_{\infty }}\end{displaymath} \ \par
\ \par
We easily obtain:\ \par

\begin{displaymath} 
\frac{1}{h}\left\vert{\left({u_{i_{0}}^{n+\frac{1}{2}}}\right) ^{q+1}-\left({u_{i_{0}-1}^{n+\frac{1}{2}}}\right) 
^{q+1}}\right\vert \leq (q+1)\left\Vert{u_{h}^{n}}\right\Vert _{\infty 
}^{q}\left\vert{v_{i_{0}}^{n+\frac{1}{2}}}\right\vert \end{displaymath} 
\ \par
and\ \par

\begin{displaymath} 
\frac{1}{h}\left\vert{\left({u_{i_{0}}^{n+\frac{1}{2}}}\right) ^{q}-\left({u_{i_{0}}^{n+\frac{1}{2}}}\right) 
^{q}}\right\vert u_{i_{0}}^{n+\frac{1}{2}}\leq \left\vert{v_{i_{0}}^{n+\frac{1}{2}}}\right\vert 
\ \left\Vert{u_{h}^{n}}\right\Vert _{\infty }^{q}\end{displaymath} \ 
\par
\ \par
and we get:\ \par

\begin{displaymath} 
\left({1-mq\Delta t_{n}\left\Vert{u_{h}^{n}}\right\Vert _{\infty }^{q}}\right) 
\ \left\Vert{v_{h}^{n+1}}\right\Vert _{\infty }\leq \end{displaymath} 
\ \par

\begin{displaymath} 
\left\Vert{v_{h}^{n}}\right\Vert _{\infty }\left({1+m(1-q^{2})\Delta t_{n}\left\Vert{u_{h}^{n}}\right\Vert 
_{\infty }^{q}+mq\Delta t_{n}\left\Vert{u_{h}^{n}}\right\Vert _{\infty 
}^{q}\frac{1+m(1-q)\Delta t_{n}\left\Vert{u_{h}^{n+1}}\right\Vert _{\infty 
}^{q}}{1-mq\Delta t_{n}\left\Vert{u_{h}^{n}}\right\Vert _{\infty }^{q}}}\right) 
\end{displaymath} \ \par
By ~(\ref{HSregime8}), we get easily there exist a positive constant 
$C\ $depending on $m,q,\ T,\ \left\Vert{u_{h}^{0}}\right\Vert _{\infty }\ $such 
that \ \par

\begin{displaymath} 
\left\Vert{v_{h}^{n+1}}\right\Vert _{\infty }\leq \left({1+C\Delta t_{n}}\right) 
\left\Vert{v_{h}^{n}}\right\Vert _{\infty }\end{displaymath} \ \par
Hence for $\displaystyle t_{n}\leq T<T_{1}$, we get $\displaystyle \displaystyle \left\Vert{\displaystyle v_{h}^{n}}\right\Vert _{\infty 
}\leq C$.\ \par
\ \par
\begin{Lemma} Under the hypotheses of Proposition ~(\ref{HSregime3}), 
we have the estimate:\ \par

\begin{equation} 
\displaystyle \left\Vert{\displaystyle v_{h}^{n}}\right\Vert _{1}\leq 
C\label{HSregime211}
\end{equation} \ \par
for $t_{n}<T<T_{1}\ $where $C$ is a constant depending on $T$ and $u_{0}$.\ 
\par
\end{Lemma}
\ \par
\textbf{Proof:  }From the properties of the semigroup operator $S$ ~\cite{mnLeRoux6}, 
we get: $\left\Vert{v_{h}^{n+\frac{1}{2}}}\right\Vert _{1}\leq \left\Vert{v_{h}^{n}}\right\Vert 
_{1}$, and by using the equations satisfied by $v_{i}^{n+\frac{1}{2}},\ -N_{-}^{n+1}\leq i\leq N_{+}^{n+1}$, 
we obtain:\ \par

\begin{displaymath} 
\left({1-mq\Delta t_{n}\left\Vert{u_{h}^{n+\frac{1}{2}}}\right\Vert _{\infty 
}^{q}}\right) \left\Vert{v_{h}^{n+1}}\right\Vert _{1}\leq \left\Vert{v_{h}^{n+\frac{1}{2}}}\right\Vert 
_{1}\left({1+m\Delta t_{n}\left\Vert{u_{h}^{n+\frac{1}{2}}}\right\Vert 
_{\infty }^{q}\left({q^{2}\left\Vert{u_{h}^{n+1}}\right\Vert _{\infty 
}+(1-q^{2})\left\Vert{u_{h}^{n+\frac{1}{2}}}\right\Vert _{\infty }}\right) 
}\right) \end{displaymath} \ \par
and we immediately deduce the estimate ~(\ref{HSregime211}).\ \par
\ \par
In this case, since $q<1,\ $ the variation of $v_{h}^{n}$ is not bounded. 
\ \par
\ \par
\section{Blow-up of the solution}
\setcounter{equation}{0}\ \par
In this part, we prove that for $q<1$ the solution blows up in finite 
time.\ \par
\ \par
\subsection{Construction of unbounded solutions}
\ \par
Define the function\ \par

\begin{displaymath} 
\theta (x)=\left\lbrace{\begin{matrix}{1-\frac{x^{2}}{a^{2}},\ \left\vert{x}\right\vert 
\leq a}\cr {0,\ \left\vert{x}\right\vert \geq a}\cr \end{matrix}}\right. 
\label{HSregime23}\end{displaymath} \ \par
\ \par
We note $\theta _{h}=\pi _{h}\theta $. If the initial condition is $u_{h}^{0}=\frac{\lambda }{T^{\frac{1}{q}}}\theta _{h}(\xi ^{0})$ 
with $\xi ^{0}=\frac{x}{T^{\frac{q-1}{2q}}}$, we prove that it is possible 
to choose $\lambda $ and $\displaystyle a$ in such a manner that $u_{h}^{n}(x)\geq \frac{\lambda }{\left({T-t_{n}}\right) ^{\frac{1}{q}}}\theta 
_{h}(\xi _{n}),\ $with\ \par
$\xi _{n}=x\left({T-t_{n})}\right) ^{\frac{1-q}{2q}}$ and then the numerical 
solution blows up in finite time.\ \par
\ \par
We denote $\displaystyle \zeta _{n}=\displaystyle \left({\displaystyle T-t_{n}}\right) ^{\displaystyle 
\frac{\displaystyle q-1}{\displaystyle 2q}}$, ${\hat u}_{h}^{n}(x)=\frac{\lambda }{\left({T-t_{n}}\right) ^{\frac{1}{q}}}\theta 
_{h}(\xi _{n})$.\ \par
The support of ${\hat u}_{h}^{n}$ is $]-a\zeta _{n},a\zeta _{n}[$; its 
lenth is increasing with the time.\ \par
If ${\hat u}_{h}^{n}\leq u_{h}^{n}$, we get ${\hat u}_{h}^{n+\frac{1}{2}}\leq u_{h}^{n+\frac{1}{2}}$.\ 
\par
The support of ${\hat u}_{h}^{n+\frac{1}{2}}$ is $\displaystyle [-\hat s_{-}^{n+1},\hat s_{+}^{n+1}]$ 
with $\hat s_{+}^{n+1}\leq s^{n+1},\ \hat s_{-}^{n+1}\leq s_{-}^{n+1}$\ 
\par
Since ${\hat u}_{h}^{n}$ is a symmetric function, it is sufficient to 
study the case $x\geq 0$,$(i\geq 0)$.\ \par
We get for $i\geq 0,$\ \par

\begin{displaymath} 
{\hat u}_{i}^{n+\frac{1}{2}}={\hat u}_{i}^{n}+\frac{\Delta t_{n}}{m}\left({\hat 
v_{i}^{n}}\right) ^{2}\end{displaymath} \ \par
for $i\ $ such that $x_{i}\leq a\zeta _{n}$ since the function ${\hat u}_{h}^{n}$ 
is decreasing for $x\geq 0$ and we have:\ \par

\begin{displaymath} 
\hat v_{i}^{n}=\frac{{\hat u}_{i}^{n}-{\hat u}_{i-1}^{n}}{h}=\frac{\lambda 
}{\left({T-t_{n}}\right) ^{\frac{1}{q}}}\frac{\theta _{i}^{n}-\theta 
_{i-1}^{n}}{h}\end{displaymath} \ \par
with $\displaystyle \theta _{i}^{n}=\theta (\xi _{i}^{n})$.\ \par
Hence, we obtain:\ \par

\begin{displaymath} 
{\hat u}_{i}^{n+\frac{1}{2}}=\frac{\lambda }{\left({T-t_{n}}\right) ^{\frac{1}{q}}}\left({\theta 
_{i}^{n}+\frac{4\lambda }{ma^{2}}\frac{\Delta t_{n}}{\left({T-t_{n}}\right) 
^{\frac{1}{q}}}\left({1-\theta _{i-\frac{1}{2}}^{n}}\right) }\right) 
\end{displaymath} \ \par
with $\theta _{i-\frac{1}{2}}^{n}=\theta (\xi ^{n}_{i-\frac{1}{2}}),\ \xi _{i-\frac{1}{2}}^{n}=\frac{1}{2}\left({\xi 
_{i}^{n}+\xi _{i-1}^{n}}\right) $.\ \par
Then ${\hat u}_{h}^{n+1}\ $will be a subsolution of ~(\ref{HSregime5}) 
if\ \par

\begin{displaymath} 
\frac{\lambda }{\left({T-t_{n+1}}\right) ^{\frac{1}{q}}}\left({1-mq\Delta 
t_{n}\left({{\hat u}_{i}^{n+\frac{1}{2}}}\right) ^{q}}\right) \theta 
_{i}^{n+1}+\frac{\lambda }{\left({T-t_{n+1}}\right) ^{\frac{1}{q}}}\frac{\Delta 
t_{n}}{h^{2}}{\hat u}_{i}^{n+\frac{1}{2}}\left({2\theta _{i}^{n+1}-\theta 
_{i-1}^{n+1}-\theta _{i+1}^{n+1}}\right) \end{displaymath} \ \par
\ \par

\begin{displaymath} 
\leq {\hat u}_{i}^{n+\frac{1}{2}}+m(1-q)\Delta t_{n}\left({{\hat u}_{i}^{n+\frac{1}{2}}}\right) 
^{q+1}\end{displaymath} \ \par
By using the equality :$\displaystyle \frac{\displaystyle 1}{\displaystyle h^{2}}\displaystyle \left({\displaystyle 
2\theta _{i}^{n+1}-\theta _{i-1}^{n+1}-\theta _{i+1}^{n+1}}\right) =\frac{\displaystyle 
2}{\displaystyle a^{2}\zeta _{n+1}^{2}}$\ \par
this inequality reduces after simplifications to\ \par

\begin{displaymath} 
\frac{\lambda }{\left({T-t_{n+1}}\right) ^{\frac{1}{q}}}\theta _{i}^{n+1}\leq 
{\hat u}_{i}^{n+\frac{1}{2}}\left({1-\frac{2\lambda \Delta t_{n}}{a^{2}\left({T-t_{n+1}}\right) 
}}\right) +\frac{\lambda }{\left({T-t_{n+1}}\right) ^{\frac{1}{q}}}mq\Delta 
t_{n}\left({{\hat u}_{i}^{n+\frac{1}{2}}}\right) ^{q}\theta _{i}^{n+1}\end{displaymath} 
\ \par

\begin{displaymath} 
+m(1-q)\Delta t_{n}\left({{\hat u}_{i}^{n+\frac{1}{2}}}\right) ^{q+1}\label{HSregime10}\end{displaymath} 
\ \par
\ \par
Noting that $\theta ^{n}_{i-\frac{1}{2}}=\theta _{i}^{n}+\frac{h}{a^{2}\zeta _{n}}\xi 
_{i-\frac{1}{2}}^{n}=\theta _{i}^{n}+\eta _{i}^{n}$ with $\displaystyle 
\displaystyle \left\vert{\displaystyle \eta _{i}^{n}}\right\vert \leq 
\frac{\displaystyle h}{\displaystyle a\zeta _{n}}$, we get:\ \par

\begin{displaymath} 
{\hat u}_{i}^{n+\frac{1}{2}}=\frac{\lambda }{\left({T-t_{n}}\right) ^{\frac{1}{q}}}\left({\theta 
_{i}^{n}\left({1-\frac{4\lambda }{ma^{2}}\frac{\Delta t_{n}}{T-t_{n}}}\right) 
+\frac{4\lambda }{ma^{2}}\frac{\Delta t_{n}}{T-t_{n}}\left({1-\eta _{i}^{n}}\right) 
}\right) \end{displaymath} \ \par
and the inequality ~(\ref{HSregime10}) becomes:\ \par

\begin{displaymath} 
\frac{\theta _{i}^{n+1}}{\left({T-t_{n+1}}\right) ^{\frac{1}{q}}}\leq 
\frac{\theta _{i}^{n}}{\left({T-t_{n}}\right) ^{\frac{1}{q}}}\left({1-\frac{4\lambda 
}{ma^{2}}\frac{\Delta t_{n}}{T-t_{n}}-\frac{2\lambda }{ma^{2}}\frac{\Delta 
t_{n}}{T-t_{n+1}}\left({1-\frac{4\lambda }{ma^{2}}\frac{\Delta t_{n}}{T-t_{n}}}\right) 
}\right) \end{displaymath} \ \par

\begin{displaymath} 
+\frac{4\lambda }{ma^{2}}\frac{\Delta t_{n}}{\left({T-t_{n+1}}\right) 
^{\frac{q+1}{q}}}\left({1-\frac{2\lambda }{ma^{2}}\frac{\Delta t_{n}}{T-t_{n+1}}}\right) 
\left({1-\eta _{i}^{n}}\right) +mq\frac{\Delta t_{n}}{\left({T-t_{n+1}}\right) 
^{\frac{1}{q}}}\left({{\hat u}_{i}^{n+\frac{1}{2}}}\right) ^{q}\theta 
_{i}^{n+1}\end{displaymath} \ \par

\begin{displaymath} 
+m(1-q)\frac{\Delta t_{n}}{\lambda }\left({{\hat u}_{i}^{n+\frac{1}{2}}}\right) 
^{q+1}\end{displaymath}  \ \par
By using the equality $\theta _{i}^{n+1}=\theta _{i}^{n}\frac{\zeta _{n}^{2}}{\zeta _{n+1}^{2}}+1-\frac{\zeta 
_{n}^{2}}{\zeta _{n+1}^{2}},$ this inequality reduces to:\ \par

\begin{displaymath} 
\theta _{i}^{n}\left({1+\frac{4\lambda }{ma^{2}}\frac{T-t_{n+1}}{T-t_{n}}+\frac{2\lambda 
}{a^{2}}\left({1-\frac{4\lambda }{ma^{2}}\frac{\Delta t_{n}}{T-t_{n}}}\right) 
}\right) \end{displaymath} \ \par

\begin{displaymath} 
\leq \frac{4\lambda }{ma^{2}}\frac{T-t_{n+1}}{T-t_{n}}\left({1-\frac{2\lambda 
}{a^{2}}\frac{\Delta t_{n}}{T-t_{n+1}}(1-\eta _{i}^{n})}\right) -\left({\frac{\zeta 
_{n+1}^{2}}{\zeta _{n}^{2}}-1}\right) \frac{T-t_{n}}{\Delta t_{n}}\end{displaymath} 
 \ \par

\begin{displaymath} 
+mq\left({T-t_{n+1}}\right) ^{1-\frac{1}{q}}\left({T-t_{n}}\right) ^{\frac{1}{q}}\left({{\hat 
u}_{i}^{n+\frac{1}{2}}}\right) ^{q}\theta _{i}^{n+1}+m(1-q)\frac{T-t_{n+1}}{\lambda 
}\left({T-t_{n}}\right) ^{\frac{1}{q}}\left({{\hat u}_{i}^{n+\frac{1}{2}}}\right) 
^{q+1}.\end{displaymath} \ \par
If we denote $\displaystyle \mu _{n}=\frac{\displaystyle 4\lambda }{\displaystyle ma^{2}}\frac{\displaystyle 
\Delta t_{n}}{\displaystyle T-t_{n}}$ since ${\hat u}_{i}^{n+\frac{1}{2}}\geq \frac{\lambda }{\left({T-t_{n}}\right) 
^{\frac{1}{q}}}\theta _{i}^{n}\left({1-\mu _{n}}\right) $, the preceding 
inequality will be satisfied if:\ \par

\begin{displaymath} 
\theta _{i}^{n}\left({1+\frac{4\lambda }{ma^{2}}+\frac{2\lambda }{a^{2}}-\mu 
_{n}\left({1+\frac{2\lambda }{a^{2}}}\right) }\right) \leq \left({\frac{4\lambda 
}{ma^{2}}-\mu _{n}\left({1+\frac{2\lambda }{a^{2}}}\right) }\right) \left({1-\eta 
_{i}^{n}}\right) -\left({\frac{\zeta _{n+1}^{2}}{\zeta _{n}^{2}}-1}\right) 
\frac{T-t_{n}}{\Delta t_{n}}\end{displaymath} \ \par

\begin{displaymath} 
+mq\lambda ^{q}\frac{\zeta _{n+1}^{2}}{\zeta _{n}^{2}}\left({\theta _{i}^{n}}\right) 
^{q}\left({1-\mu _{n}}\right) ^{q}\left({\theta _{i}^{n}\frac{\zeta _{n}^{2}}{\zeta 
_{n+1}^{2}}+1-\frac{\zeta _{n}^{2}}{\zeta _{n+1}^{2}}}\right) \end{displaymath} 
\ \par

\begin{equation} 
+m(1-q)\lambda ^{q}\frac{\displaystyle T-t_{n+1}}{\displaystyle T-t_{n}}\displaystyle 
\left({\displaystyle \theta _{i}^{n}}\right) ^{q+1}\displaystyle \left({\displaystyle 
1-\mu _{n}}\right) ^{q+1}\label{HSregime11}
\end{equation} \ \par
From the stability condition, we get:\ \par

\begin{displaymath} 
\mu _{n}\leq \frac{h}{a\zeta _{n}}\leq \frac{h}{a\zeta _{0}}=\delta \end{displaymath} 
\ \par
hence for $h$ sufficiently small, we get:$\displaystyle \frac{\displaystyle 4\lambda }{\displaystyle ma^{2}}-\mu _{n}\displaystyle 
\left({\displaystyle 1+\frac{\displaystyle 2\lambda }{\displaystyle a^{2}}}\right) 
>0,\ $ and $\left\vert{\eta _{i}^{n}}\right\vert \leq \delta $\ \par
So the inequality ~(\ref{HSregime11}) will be satisfied if :\ \par

\begin{displaymath} 
\theta _{i}^{n}\left({1+\frac{4\lambda }{ma^{2}}+\frac{2\lambda }{a^{2}}-\mu 
_{n}\left({1+\frac{2\lambda }{a^{2}}}\right) }\right) \leq \left({\frac{4\lambda 
}{ma^{2}}-\mu _{n}\left({1+\frac{2\lambda }{a^{2}}}\right) }\right) \left({1-\delta 
}\right) -\left({\frac{\zeta _{n+1}^{2}}{\zeta _{n}^{2}}-1}\right) \frac{T-t_{n}}{\Delta 
t_{n}}\end{displaymath} \ \par

\begin{displaymath} 
+mq\lambda ^{q}\left({1-\mu _{n}}\right) ^{q}\left({\theta _{i}^{n}}\right) 
^{q}\left({\theta _{i}^{n}+\frac{\zeta _{n+1}^{2}}{\zeta _{n}^{2}}-1}\right) 
+m(1-q)\lambda ^{q}\left({1-\mu _{n}}\right) ^{q+1}\left({\theta _{i}^{n}}\right) 
^{q+1}\frac{T-t_{n+1}}{T-t_{n}}\end{displaymath} \ \par
Since $\theta _{i}^{n}\in (0,1),$ we introduce the function $\Phi _{n}(y)$ 
defined on $(0,1)$ by :\ \par

\begin{displaymath} 
\Phi _{n}(y)=m\lambda ^{q}\left({1-\mu _{n}}\right) ^{q}y^{q+1}\left({q+(1-q)(1-\mu 
_{n})\frac{T-t_{n+1}}{T-t_{n}}}\right) \end{displaymath} \ \par

\begin{displaymath} 
+\left({\frac{4\lambda }{ma^{2}}-\mu _{n}\left({1+\frac{2\lambda }{a^{2}}}\right) 
}\right) (1-\delta )-\left({\frac{\zeta _{n+1}^{2}}{\zeta _{n}^{2}}-1}\right) 
\frac{T-t_{n}}{\Delta t_{n}}\end{displaymath} \ \par

\begin{displaymath} 
-y\left({1+\frac{4\lambda }{ma^{2}}+\frac{2\lambda }{a^{2}}-\mu _{n}\left({1+\frac{2\lambda 
}{a^{2}}}\right) }\right) \end{displaymath} \ \par
or $\displaystyle \Phi _{n}(y)=A\lambda ^{q}y^{q+1}+C-By$\ \par
with 
\begin{displaymath} 
A=m\left({1-\mu _{n}}\right) ^{q}\left({q+(1-q)(1-\mu _{n})\frac{T-t_{n+1}}{T-t_{n}}}\right) 
\end{displaymath} \ \par

\begin{displaymath} 
B=1+\frac{4\lambda }{ma^{2}}+\frac{2\lambda }{a^{2}}-\mu _{n}\left({1+\frac{2\lambda 
}{a^{2}}}\right) \end{displaymath} \ \par

\begin{displaymath} 
C=\left({\frac{4\lambda }{ma^{2}}-\mu _{n}\left({1+\frac{2\lambda }{a^{2}}}\right) 
}\right) (1-\delta )-\left({\frac{\zeta _{n+1}^{2}}{\zeta _{n}^{2}}-1}\right) 
\frac{T-t_{n}}{\Delta t_{n}}\end{displaymath} \ \par
A sufficient condition to satisfy ~(\ref{HSregime11}) is: $\displaystyle 
\Phi _{n}(y)\geq 0,\ y\in (0,1).$\ \par
We have $\Phi _{n}(0)=C$,\ \par
hence we get\ \par

\begin{displaymath} 
\Phi _{n}(0)\geq \left({\frac{4\lambda }{ma^{2}}-\delta \left({1+\frac{2\lambda 
}{a^{2}}}\right) }\right) (1-\delta )-\left({\frac{\zeta _{n+1}^{2}}{\zeta 
_{n}^{2}}-1}\right) \frac{T-t_{n}}{\Delta t_{n}}\end{displaymath} \ \par
But since $0<q<1$, we get:\ \par

\begin{displaymath} 
0\leq \frac{\zeta _{n+1}^{2}}{\zeta _{n}^{2}}-1\leq \frac{1-q}{q}\frac{\Delta 
t_{n}}{T-t_{n+1}}\frac{\zeta _{n}^{2}}{\zeta _{n+1}^{2}}\end{displaymath} 
\ \par
and 
\begin{displaymath} 
C\geq \left({\frac{4\lambda }{ma^{2}}-\delta \left({1+\frac{2\lambda }{a^{2}}}\right) 
}\right) (1-\delta )-\frac{1-q}{q}.\end{displaymath} \ \par
So, if the quantity $\frac{\lambda }{a^{2}}$ is such that $\displaystyle 
\frac{\displaystyle 4\lambda }{\displaystyle ma^{2}}>\frac{\displaystyle 
1-q}{\displaystyle q},\ $if $h\ $is sufficiently small, we get $C>0$.\ 
\par
\ \par
Let us define $y_{0}=\frac{C}{B}$, $y_{0}\in (0,1)$ and if $y\leq y_{0}$, 
we obtain $\displaystyle \Phi _{n}(y)\geq 0$; if $y_{0}\leq y\leq 1,\ $we 
obtain:\ \par

\begin{displaymath} 
\Phi _{n}(y)\geq A\lambda ^{q}y_{0}^{q+1}-B(1-y_{0})\end{displaymath} 
.\ \par
This quantity will be positive if : $\lambda ^{q}\geq \frac{B-C}{Ay_{0}^{q+1}}$.\ 
\par
Further, we have:\ \par

\begin{displaymath} 
A\geq m(1-\delta )^{q}\left({q+(1-q)(1-\delta )\left({1-\frac{\Delta t_{n}}{T-t_{n}}}\right) 
}\right) \end{displaymath} \ \par
The solution at the time level $t_{n+1}$ exists if $mq\Delta t_{n}\left\Vert{{\hat u}_{h}^{n}}\right\Vert _{\infty }<1$, 
that is $\displaystyle \frac{\displaystyle \Delta t_{n}}{\displaystyle T-t_{n}}<\frac{\displaystyle 
1}{\displaystyle mq\lambda ^{q}}$\ \par
and if $\lambda \geq \lambda _{0}>\frac{1}{\left({mq}\right) ^{\frac{1}{q}}}$, 
we get: $\displaystyle A\geq mq(1-\delta )^{q}$\ \par
and $\Phi _{n}(y)$ will be positive if:\ \par

\begin{equation} 
\lambda ^{q}\geq \frac{\displaystyle 1}{\displaystyle mq(1-\delta )^{q}}\displaystyle 
\left({\displaystyle \frac{\displaystyle 1}{\displaystyle q}+\frac{\displaystyle 
2\lambda }{\displaystyle a^{2}}+\frac{\displaystyle 4\lambda }{\displaystyle 
ma^{2}}\delta }\right) \displaystyle \left({\displaystyle \frac{\displaystyle 
1+\frac{\displaystyle 4\lambda }{\displaystyle ma^{2}}+\frac{\displaystyle 
2\lambda }{\displaystyle a^{2}}}{\displaystyle \frac{\displaystyle 4\lambda 
}{\displaystyle ma^{2}}(1-\delta )-\delta \displaystyle \left({\displaystyle 
1+\frac{\displaystyle 2\lambda }{\displaystyle ma^{2}}}\right) -\frac{\displaystyle 
1-q}{\displaystyle q}}}\right) ^{q+1}\label{HSregime12}
\end{equation} \ \par
The second member is a function of $\frac{\lambda }{a^{2}}$, hence if 
$\frac{\lambda }{a^{2}}$ is fixed such that $\displaystyle \frac{\displaystyle 4\lambda }{\displaystyle ma^{2}}>\frac{\displaystyle 
1-q}{\displaystyle q}$,\ \par
the inequality ~(\ref{HSregime12}) will be satisfied if $\lambda $ is 
large enough and $h$ sufficiently small.\ \par
Hence if the initial condition satisfies $u_{0h}(x)\geq \frac{\lambda }{T^{\frac{1}{q}}}\theta _{h}(\xi _{0})$, 
$\lambda $ and $a$ satisfying ~(\ref{HSregime12}), the solution blows 
up in finite time.\ \par
\ \par
\subsection{Blow-up of the solution}
\ \par
\begin{Theorem} Let $0<q<1.$The solution of problem ~(\ref{HSregime1}) 
blows up in finite time.\ \par
\end{Theorem}
\ \par
\textbf{Proof:}Since $u_{0}(x)\not\equiv 0$, there exists $\rho >0,\ \epsilon >0,\ x_{0}$ 
such that $u_{0}(x)\geq \epsilon >0$ for $x\ $such that\ \par
$\left\vert{x-x_{0}}\right\vert <\rho .\ $So we can choose $T$ large 
enough such that $\frac{\lambda }{T^{\frac{1}{q}}}<\epsilon $ and $\frac{a}{T^{\frac{1-q}{2q}}}<\rho $.\ 
\par
We have $u_{0}(x)\geq \frac{\lambda }{T^{\frac{1}{q}}}\theta _{h}\left({\frac{x-x_{0}}{\zeta 
_{0}}}\right) $ and then the solution blows up at time $T_{0}\leq T^{*}$ 
where\ \par
$T^{*}={\mathrm{Max}}\left({\left({\frac{\lambda }{\epsilon }}\right) ^{q},\left({\frac{a}{\rho 
}}\right) ^{\frac{2q}{1-q}}}\right) $.\ \par
\ \par
\ \par
In Fig1, we present the evolution of an initial condition $u_{0}\ $ for 
$m=1,\ p=1.5$. The solution blows up in finite time.\ \par
\ \par

\begin{figure}[h]
\begin{center}
\rotatebox{0}{\resizebox{10cm}{!}{\includegraphics{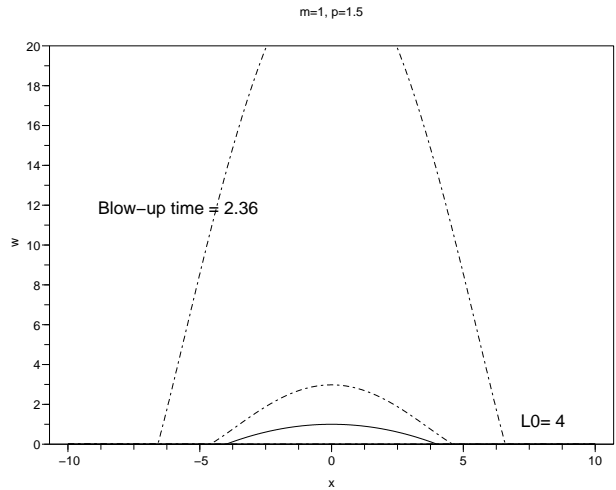}}}
\caption{ blow-up for m=1, p=1,5}
\label{HSregime100}
\end{center}
\end{figure}
\ \par
\ \par

\bibliographystyle{C:/TexLive/texmf/bibtex/bst/base/plain}

\end{document}